\documentclass[a4paper,12pt]{amsart}

\nonstopmode \numberwithin{equation}{section}
\usepackage[english]{babel}
\usepackage{amssymb}
\usepackage{amsmath}

\newtheorem{theorem}{Theorem}[section]
\newtheorem{lemma}[theorem]{Lemma}

\numberwithin{equation}{section}

\newcommand{\C}{{\mathbb C}}
\newcommand{\class}{{\mathcal C}}
\newcommand{\sclass}{{\widehat{\mathcal C}}}
\newcommand{\Ch}{{\widehat{C}}}
\newcommand{\D}{{\mathbb D}}

\newcommand{\uhp}{{\mathbb H}}
\newcommand{\J}{{\mathcal J}}

\newcommand{\R}{{\mathbb R}}
\newcommand{\Rh}{{\widehat{\mathbb R}}}

\newcommand{\Cd}{{\widetilde C}}

\newcommand{\Ce}{{\widehat C}}

\newcommand{\sphere}{{\widehat{\mathbb C}}}

\renewcommand{\Im}{\,{\operatorname{Im}\,}}

\renewcommand{\Re}{{\operatorname{Re}\,}}

\newcommand{\inv}{^{-1}}

\renewcommand{\arg}{\,{\operatorname{arg}\,}}

\begin{document}
\bibliographystyle{amsplain}

\title{Spherical convexity and hyperbolic metric}
\dedicatory{Dedicated to Professor David Minda \\
on the occasion of his retirement}


\author{Toshiyuki Sugawa}
\address{Graduate School of Information Sciences, Tohoku University \\
Sendai 980-8579, Japan}
\email{sugawa@math.is.tohoku.ac.jp}

\maketitle


\begin{abstract} 
Let $\Omega$ be a domain in $\C$ with hyperbolic metric
$\lambda_\Omega(z)|dz|$ with Gaussian curvature $-4.$
Mej\'\i a and Minda proved in their 1990 paper that $\Omega$ is
(Euclidean) convex if and only if $d(z,\partial\Omega)\lambda_\Omega(z)\ge1/2$
for $z\in\Omega,$ where $d(z,\partial\Omega)$ 
denotes the Euclidean distance
from $z$ to the boundary $\partial\Omega.$
In the present note, we will provide similar characterizations of
spherically convex domains in terms of the spherical density of
the hyperbolic metric.
\end{abstract}

\section{Introduction}

Let $\Omega$ be a domain in the Riemann sphere $\sphere=\C\cup\{\infty\}.$
The Uniformization Theorem implies that there exists a holomorphic
universal covering projection $p_\Omega$
of the unit disk $\D=\{\zeta\in\C: |\zeta|<1\}$
onto $\Omega$ whenever $\partial\Omega$ contains at least three points.
Here and hereafter, the boundary of a domain is taken in the Riemann sphere.
For instance, $\partial\C=\{\infty\}.$
Since the Poincar\'e metric $|d\zeta|/(1-|\zeta|^2)$ is invariant
under analytic automorphisms of $\D,$ the quantity
$$
\lambda_\Omega(z)=\frac{1}{(1-|\zeta|^2)|p_\Omega'(\zeta)|},
\quad z=p_\Omega(\zeta),
$$
does not depend on the particular choice of the covering projection
$p_\Omega$ and the pre-image $\zeta\in p_\Omega\inv(z).$
It is well known that $\lambda_\Omega(z)|dz|$ gives $\Omega$
a complete metric with constant Gaussian curvature $-4.$
The metric $\lambda_\Omega(z)|dz|$ is called {\it hyperbolic}
and $\Omega$ is also called hyperbolic if it admits such a metric.
The Schwarz lemma yields the monotonicity of the hyperbolic metric;
namely, $\lambda_\Omega(z)\le \lambda_{\Omega'}(z)$ whenever
$\Omega'\subset\Omega.$
For instance, $\lambda_D(z)=r/(r^2-|z-a|^2)$ for the disk
$D=\{z: |z-a|<r\}.$
Let
$$
d(z,\partial\Omega)=\inf_{a\in\partial\Omega}|z-a|,\quad z\in\Omega,
$$
Then $\Delta=\{w: |w-z|<d(z,\partial\Omega)\}$ 
is the maximal disk centered at $z$ and contained in $\Omega.$
Thus the monotonicity of the hyperbolic metric implies that
$\lambda_\Omega(z)\le\lambda_\Delta(z)=1/d(z,\partial\Omega)$ and, equivalently,
$d(z,\partial\Omega)\lambda_\Omega(z)\le 1.$
On the other hand, the Koebe one-quarter theorem implies
$d(z,\partial\Omega)\lambda_\Omega(z)\ge 1/4$ for a simply connected domain 
$\Omega.$
In general, however, such a positive lower bound does not necessarily exist.
Consider the quantity
$$
C(\Omega)=\inf_{z\in\Omega}d(z,\partial\Omega)\lambda_\Omega(z).
$$
Beardon and Pommerenke \cite{BP78} first observed that
domains $\Omega$ with $C(\Omega)>0$ are characterized by uniform perfectness
of its boundary.
Here, a closed set $A$ in $\sphere$ consisting at least two points
is said to be {\it uniformly perfect} if there exists a constant $0<c<1$
such that for any point $a\in A\setminus\{\infty\}$ and $0<r<d(A),$ 
one can find a point $b\in A$ with $cr\le |a-b|\le r.$
Here, $d(A)$ denotes the Euclidean diameter of $A.$
We define $d(A)$ to be $+\infty$ when $\infty\in A.$
It is easy to see that a uniformly perfect set is perfect
and hence contains uncountably many points.
It is also known that
uniform perfectness is invariant under M\"obius transformations.
Many characterizing properties of uniformly perfect sets are known.
The reader can consult recent monographs \cite{AW:SP} or \cite{KL:hg}
for hyperbolic metric and uniformly perfect sets.
The following is due to Beardon and Pommerenke \cite{BP78}.

\begin{theorem}[Beardon-Pommerenke]\label{thm:BP}
Let $\Omega$ be a hyperbolic domain in $\C.$
Then $\partial\Omega$ is uniformly perfect if and only if
$C(\Omega)>0.$
\end{theorem}

Concerning the value of $C(\Omega),$ the following beautiful results should be
recalled.

\begin{theorem}[Harmelin-Minda-Mej\'\i a]\label{thm:M}
Let $\Omega$ be a hyperbolic domain in $\C.$
Then $C(\Omega)\le 1/2$ and equality holds precisely when
$\Omega$ is convex.
\end{theorem}

The inequality $C(\Omega)\ge 1/2$ for a convex domain $\Omega$
follows easily from the classical covering theorem for convex functions
(see \cite[Theorem 2.15]{Duren:univ} for instance).
On the other hand, the converse looks subtle because the domain is not
assumed to be simply connected.
Mej\'\i a and Minda \cite[Theorem 2]{MM90} gave a proof of the converse. 
(They remark that the first proof was done by Hilditch.)
The general inequality $C(\Omega)\le 1/2$ is due to Harmelin and Minda
\cite{HM92}.

We note that Theorem \ref{thm:BP} is valid only for domains in $\C$
(see, for instance, \cite{Sugawa13}).
To include domains in $\sphere,$ it is natural to consider
spherical analogs.
When we deal with subsets of $\sphere,$ it is convenient to use
the chordal distance
\begin{equation}\label{eq:dist}
\sigma(z,w)=\frac{|z-w|}{\sqrt{(1+|z|^2)(1+|w|^2)}}
=\frac{\tau(z,w)}{\sqrt{1+\tau(z,w)^2}},
\end{equation}
where
$$
\tau(z,w)=\left|\frac{z-w}{1+z\bar w}\right|.
$$
When one of $z,w$ is $\infty,$ we should take a suitable limit so that
$\sigma(z,\infty)=1/\sqrt{1+|z|^2}$ for instance.
Note here that $\sigma(z,w)\le\tau(z,w).$
Though $\tau$ does not satisfy the triangle inequality,
it is often a good substitute of the Euclidean distance 
$d(z,w)=|z-w|$
for spherical geometry (cf.~\cite{Minda85a} or \cite{Sugawa13}).
We observe that $\sigma$ and $\tau$ both generate
the spherical metric $\sigma(z)|dz|=|dz|/(1+|z|^2):$
$$
\lim_{w\to z}\frac{\sigma(z,w)}{|z-w|}
=\lim_{w\to z}\frac{\tau(z,w)}{|z-w|}
=\frac1{1+|z|^2}=\sigma(z).
$$
The geodesic distance $\theta(z,w)$ induced by the Riemannian metric 
$\sigma(z)|dz|$ is called the spherical distance and explicitly given by
$$
\theta(z,w)=\arcsin\sigma(z,w)=\arctan\tau(z,w),
$$
which is the length of shorter great circle $\gamma$ joining $z$ and $w.$
Since $\sin x<x<\tan x$ for $0<x<\pi/2,$ we obtain
\begin{equation}\label{eq:3}
\sigma(z,w)\le \theta(z,w)\le\tau(z,w)
\end{equation}
for $z,w\in\sphere.$
Geometrically, via the stereographic projection onto the sphere,
$\theta(z,w)$ is half the angle subtended by $\gamma$ at the center
of the sphere.
We can measure the distance from a point $z\in\Omega$ to the boundary
$\partial\Omega$ by $\sigma, \theta$ or $\tau.$
We set
$$
\chi(z,\partial\Omega)=\min_{a\in\partial\Omega}\chi(z,a)
$$
for $z\in\Omega,$ where $\chi$ represents one of the quantities
$\sigma, \theta$ and $\tau.$

Minda \cite{Minda85a} considered the spherical density of the
hyperbolic metric
$$
\mu_\Omega(z)=\frac{\lambda_\Omega(z)|dz|}{\sigma(z)|dz|}
=(1+|z|^2)\lambda_\Omega(z),
$$
which is invariant under the spherical isometries.
As spherical analogs of $C(\Omega),$ we consider the quantity
$$
\Ch_\chi(\Omega)=\inf_{z\in\Omega}\chi(z,\partial\Omega)\mu_\Omega(z)
$$
for $\chi=\sigma, \theta, \tau.$
By \eqref{eq:3}, we have
\begin{equation}\label{eq:Ch}
\Ch_\sigma(\Omega)\le \Ch_\theta(\Omega)\le \Ch_\tau(\Omega).
\end{equation}
Note also that these quantities are invariant under spherical isometries;
in other words, $\Ch_\chi(T(\Omega))=\Ch_\chi(\Omega)$ for 
$\chi=\sigma, \theta, \tau$ and $T(z)=e^{it}(z-a)/(1+\bar a z)$
or $T(z)=e^{it}/z$ for $t\in\R$ and $a\in\C.$
The author proved the following results in \cite{Sugawa13}
as a spherical analog of Theorem \ref{thm:BP}.
(Note that the quantities $\Ch_\sigma(\Omega)$ and $\Ch_\tau(\Omega)$
are denoted by $\Cd(\Omega)$ and $\Ce(\Omega),$ respectively, 
in \cite{Sugawa13}.
The quantity $\Ch_\theta(\Omega)$ is not defined in \cite{Sugawa13}
but the following assertion is clear from \eqref{eq:Ch}.)

\begin{theorem}[$\text{\cite[Theorem 1.4, Corollary 1.5]{Sugawa13}}$]%
\label{thm:S}
Let $\Omega$ be a hyperbolic domain in $\sphere.$
Then $\Ch_\tau(\Omega)\le 1/2.$
Moreover, the following conditions are equivalent:
\begin{enumerate}
\item
$\partial\Omega$ is uniformly perfect.
\item
$\Ch_\sigma(\Omega)>0.$
\item
$\Ch_\theta(\Omega)>0.$
\item
$\Ch_\tau(\Omega)>0.$
\end{enumerate}
\end{theorem}

We next recall the notion of spherical convexity.
A domain $\Omega$ in $\sphere$ is called {\it spherically convex}
if every pair of points in $\Omega$ can be joined by a spherical
geodesic (a subarc of a great circle) which lies entirely in $\Omega.$
One can easily observe that a spherically convex domain in $\sphere$
which contains a pair of antipodal points must be the whole sphere.
To avoid this trivial case, in what follows, we will always assume
that a spherically convex domain is a proper subdomain of $\sphere.$
Thus for a spherically convex domain $\Omega,$ one has $-1/\bar z\notin\Omega$
whenever $z\in\Omega.$
Many characteristic properties of spherically convex domains are known.
See \cite{KM01} for instance.

We are in a position to present a spherical analog of Theorem \ref{thm:M}.

\begin{theorem}\label{thm:main}
Let $\Omega$ be a hyperbolic domain in $\sphere.$
Then the following conditions are equivalent:
\begin{enumerate}
\item
$\Omega$ is spherically convex.
\item
$\Ch_\sigma(\Omega)=1/2.$
\item
$\Ch_\theta(\Omega)=1/2.$
\item
$\Ch_\tau(\Omega)=1/2.$
\end{enumerate}
\end{theorem}

In the next section, we give another characterization of
spherically convex domains, which is analogous to that of
Euclidean convex domains due to Keogh.
The final section will be devoted to the proof of Theorem \ref{thm:main}.

\section{A characterization of spherically convex domains}

The most difficult part of Theorem \ref{thm:M} is the assertion that
the inequality $C(\Omega)\ge1/2$ implies convexity of $\Omega.$
Mej\'\i a and Minda proved it by making clever use of Keogh's theorem
\cite[Theorem 1]{Keogh76}.
The proof of Keogh's theorem is based on the following lemma.
Let $\class$ be the class of domains $G$ of the form $\Delta_1\setminus 
\overline{\Delta_2}$ for open disks $\Delta_1$ and $\Delta_2$ in $\C$
whose boundaries intersect orthogonally.
The midpoint of the concave boundary arc $\partial\Delta_2\cap\Delta_1$
of $G$ will be denoted by $m(G).$

\begin{lemma}[Keogh \cite{Keogh76}]\label{lem:K}
Let $\Omega$ be a domain in $\C.$
Then, $\Omega$ is not convex if and only if there is a domain $G\in\class$
such that $G\subset\Omega$ with $m(G)\in\partial\Omega.$
\end{lemma}

Though it is hard to formulate a spherical analog of Keogh's theorem,
one can obtain that of the above lemma.
Note that the spherical disk $\{z: \tau(z,a)<\rho\}=\{z: \sigma(z,a)
<\rho/\sqrt{1+\rho^2}\}$ is spherically convex precisely if
$\rho\le 1.$
Let $\sclass$ be the class of domains $G$ of the form
$\Delta_1\setminus\overline{\Delta_2}$ for spherical disks
$\Delta_1,\Delta_2$ of $\tau$-radius less than 1,
whose boundaries intersect orthogonally.
Note that (the closure of) a spherical disk of $\tau$-radius $<1$
is strictly spherically convex.
Note also that $\tau$-radius $<1$ if and only if $\theta$-radius $<\pi/4.$
Let $\hat m(G)$ denote the spherical midpoint of the spherically concave
boundary arc $\partial\Delta_2\cap\Delta_1.$
Then, we have a spherical analog of Lemma \ref{lem:K}.

\begin{lemma}\label{lem:sc}
Let $\Omega$ be a domain in $\sphere.$
Then, $\Omega$ is not spherically convex if and only if 
there is a domain $G\in\sclass$
such that $G\subset\Omega$ with $\hat m(G)\in\partial\Omega.$
\end{lemma}

\proof
First we suppose that $G\subset\Omega$ and that $a=\hat m(G)\in\partial\Omega$
for some $G=\Delta_1\setminus\overline{\Delta_2}\in\sclass.$
Let $C$ be the great circle tangent to $\partial\Delta_2$ at $a.$
Then $C\cap \overline{\Delta_2}=\{a\}.$
Therefore, we can take a small enough open segment $\gamma$ of $C$
so that $a\in \gamma\subset\Delta_1.$
Since $\gamma\setminus\{a\}\subset G\subset\Omega,$
the domain $\Omega$ is not spherically convex.

Suppose next that $\Omega$ is not spherically convex.
Then, there is a pair of points $z_1, z_2\in \Omega$ with
$\tau(z_1,z_2)<1$ such that the shorter spherical geodesic joining
$z_1$ and $z_2$ contains a point $b$ of the complement 
$A=\sphere\setminus\Omega.$
We need to find a $G\in\sclass$ such that $G\subset\Omega$ and that
$\hat m(G)\in\partial\Omega.$
We may assume that $z_1, z_2\in\R$ and that $z_1<b<z_2.$
Take a polygonal Jordan arc $\gamma$ joining $z_1$ and $z_2$
in $\Omega$ in such a way that
$\gamma$ intersects $\Rh$ at finitely many points transversally.
We may further assume that $\gamma$ consists of only horizontal
or vertical segments.
Then the set $\sphere\setminus(\gamma\cup\Rh)$ consists of finitely many
connected components.
Let $V$ denote the set of those components.
Then $V$ consists of at least three elements.
Note that each $D\in V$ is a Jordan domain whose boundary consists of
finitely many pieces of $\gamma$ and $\Rh.$
We choose $\gamma$ so that the number $\sharp(\gamma\cap\Rh)$ of intersection
points of $\gamma$ and $\Rh$ is smallest possible.

We now construct an undirected graph for which $V$ is the set of vertices.
Let $\J$ be the collection of open intervals which are connected components
of the set $\Rh\setminus\gamma:$ the circle minus the intersection points of
$\gamma$ and $\Rh.$
We will connect distinct $D, D'\in V$ with one edge if $D$ and $D'$ share 
an interval in $\J.$
In other words, the set $E$ of edges consists of un-orderd pairs
$\{D, D'\}$ of distinct elements in $V$ such that 
$I\subset \overline{D}\cap\overline{D'}$ for some $I\in\J.$
Since $\sphere\setminus\gamma$ is connected,
$\Gamma=(V,E)$ is a connected graph.
We now show that the graph $\Gamma$ is a tree; namely, there is no
cycle in $\Gamma.$
Suppose, to the contrary, that there is a finite sequence
of distinct elements $D_1, D_2, \dots, D_n$ in $V$ with $n\ge3$ such that
$\{D_{j-1},D_j\}\in E$ for $j=1,2,\dots, n,$ where we set
$D_0=D_n$ for convenience.
By definition, there is an $I_j\in\J$ with
$I_j\subset\overline{D_{j-1}}\cap\overline{D_j}$ for each $j=1,2,\dots,n.$
Take a point $c_j$ from $I_j$ for each $j.$
We draw a closed Jordan arc $\delta_j$ joining $c_j$ and $c_{j+1}$
in such a way that $\delta_j\setminus\{c_j,c_{j+1}\}\subset D_j$ for each $j.$
Since $D_j$'s are mutually disjoint, the union $\delta=\delta_1+\delta_2+
\dots+\delta_n$ is a Jordan curve with $\gamma\cap\delta=\emptyset.$
By the Jordan curve theorem, $\sphere\setminus\delta$ consists of
two connected components, say, $W_1$ and $W_2.$
The arc $\gamma$ is contained entirely in one of the components because of its
connectedness.
We may assume that $\gamma\subset W_1.$
On the other hand, the intervals $I_j$ and $I_{j+1}$ are separated by
intersection points of $\gamma$ and $\Rh$ in $\partial D_j$
so that $\partial D_j\cap W_2 \ne\emptyset,$ 
which is a contradiction.
Thus we have shown that there is no cycle in the graph $\Gamma.$

We recall that a vertex of degree $1$ (namely, a vertex with only one edge)
is called a {\it leaf}.
Since every tree containing at least three vertices
has at least two leaves as the ends of one of the longest paths 
(see, for example, \cite[\S 1.5]{Diestel:gt}),
there are two components $D, D'\in V$ such that the boundaries of $D$ and $D'$
contain unique intervals $I$ and $I'$ in $\J,$ respectively.
Since the total length of the great circle $\Rh$ is $\pi,$
at least one of $I$ and $I'$ has length less than $\pi/2.$
We may assume that $I$ is of the form $(x_1,x_2)$ with 
$-\infty<x_1<x_2<+\infty$ and $\theta(x_1,x_2)<\pi/2.$
The other piece $\gamma_1:=\partial D\setminus I$ of the boundary
of $D$ is a closed subarc of $\gamma.$
Since $\gamma_1$ does not intersect $\Rh$ except for its endpoints,
$D$ is contained in the upper half-plane or the lower half-plane.
We will further assume that $D\subset\{z: \Im z>0\}$ because there is no 
essential difference.
We now show that $(D\cup I)\cap A\ne\emptyset.$
Suppose that $(D\cup I)\cap A$ is an empty set.
In the case when $\{x_1,x_2\}\cap\{z_1,z_2\}=\emptyset,$ we modify $\gamma$ to
$\gamma'$ by replacing $[x_1,x_1-i\eta]\cup\gamma_1\cup[x_2,x_2-i\eta]$
by the segment $[x_1-i\eta,x_2-i\eta]$ for small enough $\eta>0.$
Then, we get $\sharp(\gamma'\cap\Rh)=\sharp(\gamma\cap\Rh)-2,$
which contradicts minimality of $\sharp(\gamma\cap\Rh).$
Next consider the case when $x_1=z_1.$
Since $x_2\ne z_2$ by the assumption $(z_1,z_2)\cap A\ne\emptyset,$
we have $x_1<x_2<z_2.$
Thus, we can modify $\gamma$ to $\gamma'$ by replacing 
$\gamma_1\cup[x_2,x_2-i\eta]$ by $[x_1,x_1-i\eta]\cup[x_1-i\eta,x_2-i\eta]$
for small enough $\eta>0.$
Then, $\sharp(\gamma'\cap\Rh)=\sharp(\gamma\cap\Rh)-1,$
which contradicts minimality of $\sharp(\gamma\cap\Rh).$
The other three cases $x_1=z_2, x_2=z_1, x_2=z_2$ can be treated similarly.
We have thus shown the claim $(D\cup I)\cap A\ne\emptyset.$

We show existence of a desired $G\in\sclass$ under the assumption
that $D\cap A\ne\emptyset.$
To this end, we consider the function
$$
u(z)=\arg\frac{x_2-z}{z-x_1}
$$
which is harmonic on a neighbourhood of $D\cup I$ such that $u=0$ on $I.$
Note that $0\le u<\pi$ on $D\cup I.$
Choose a point $a\in D\cap A$ such that $u(z)$ takes its maximum at $z=a$
on the set $D\cap A.$
Let $C$ denote the circle passing through the three points
$a, x_1$ and $x_2$ and let $\Delta_2$ be the interior of $C.$
Since the exterior of $C$ contains the part $\Rh\setminus[x_1,x_2]$
of the great circle $\Rh$ with spherical length
$\pi-\theta(x_1,x_2)>\pi/2,$ the disk $\Delta_2$ is spherically convex.
We now choose a small enough disk $\Delta_1$ contained in the upper
half-plane in such a way that $\partial\Delta_1$ is orthogonal to
$C$ and that $a$ is the spherical midpoint of the arc $C\cap\Delta_1.$
Then $G=\Delta_1\setminus\overline{\Delta_2}\in\sclass$ satisfies
$G\subset\Omega$ and $\hat m(G)=a\in\partial\Omega.$

Next we assume that $D\cap A=\emptyset$ but $I\cap A\ne\emptyset.$
We take a small enough $\eta>0$ so that $\theta(x_1',x_2')<\pi/2$
and $[x_1,x_1']\cup [x_2,x_2']\subset\Omega,$ where
$x_1'=x_1-i\eta$ and $x_2'=x_2-i\eta.$
Note that the great circle passing through $x_1'$ and $x_2'$
is a Euclidean circle whose center lies in the half-plane $\Im z>-\eta.$
We consider the harmonic function
$$
u(z)=\arg\frac{x_2'-z}{z-x_1'}
$$
and take a maximum point $z=a$ of $u(z)$ on $I.$
Let $C$ denote the circle passing through the three points
$a, x_1'$ and $x_2'.$
Then the (Euclidean) center lies in the half-plane $\Im z<-\eta.$
Therefore its exterior contains the longer great circle
joining $x_1'$ and $x_2',$ which implies that the interior of $C$
is strictly spherically convex.
We can now construct $G$ in the same manner as before.
The proof of the lemma is now complete.
\qed

\section{Proof of Theorem \ref{thm:main}}

For the proof, we recall a result due to Minda.

\begin{lemma}[Minda $\text{\cite[Theorem 1]{Minda86}}$]\label{lem:M}
Let $\Omega\ne\sphere$ be a spherically convex domain in $\sphere$ 
and $z\in\Omega.$
Then
$$
\mu_\Omega(z)\ge\frac{1+\tau(z,\partial\Omega)^2}{2\tau(z,\partial\Omega)},
$$
where equality holds if and only if $\Omega$ is a hemisphere.
\end{lemma}

We now prove Theorem \ref{thm:main}.
Let $\Omega\ne\sphere$ be a spherically convex domain.
By Lemma \ref{lem:M} and \eqref{eq:dist}, we obtain
$$
\sigma(z,\partial\Omega)\mu_\Omega(z)
\ge\sigma(z,\partial\Omega)\cdot
\frac{1+\tau(z,\partial\Omega)^2}{2\tau(z,\partial\Omega)}
=\frac{\tau(z,\partial\Omega)}{2\sigma(z,\partial\Omega)}
\ge\frac12.
$$
Since $\Ch_\sigma(\Omega)\le \Ch_\theta(\Omega)\le \Ch_\tau(\Omega)
\le 1/2$ by \eqref{eq:Ch} and Theorem \ref{thm:S},
we see that $(i)\Rightarrow (ii)\Rightarrow (iii)\Rightarrow(iv).$
(In \cite{Sugawa13}, only $\Ch_\tau(\Omega)=1/2$ was shown
for a spherically convex domain $\Omega.$)

Finally, we show the implication $(iv)\Rightarrow (i).$
To this end, it suffices to show that $\Ch_\tau(\Omega)<1/2$
if a hyperbolic domain $\Omega$ is not spherically convex.
Suppose that $\Omega$ is not spherically convex.
By Lemma \ref{lem:sc}, there is a domain $G=
\Delta_1\setminus\overline{\Delta_2}\subset\Omega$
with $\hat m(G)\in\partial\Omega$
for strictly spherically convex disks $\Delta_1,\Delta_2$
whose boundaries are orthogonal to each other.
By using a spherical isometry, we may assume that
$\Delta_2=\{z: \Re z<-c\}$ for some $c\in\R$ and 
$\Delta_1=\{z: |z+c|<r\}$ for some $r>0.$
Here, we have $c>0$ because of strict spherical convexity of $\Delta_2.$
We also observe that the antipodal point $-1/(r-c)$ of $r-c\in\partial\Delta_1$
should lie in the exterior of $\Delta_1.$
This is equivalent to the condition
\begin{equation}\label{eq:cond}
r^2<c^2+1.
\end{equation}
We now compute the hyperbolic density $\lambda_G(z).$
Put $\alpha=-c+ir.$ 
Then, the M\"obius transformation $w=L(z)=(z-\alpha)/(z-\bar\alpha)$
maps $G$ onto the first quadrant $\{w\in\C: \Re w>0, \Im w>0\}$ of the
complex plane.
Hence, $\zeta=h(z)=w^2$ is a conformal mapping of $G$ onto the upper
half-plane $\uhp.$
Since $\lambda_\uhp(\zeta)=1/(2\Im\zeta),$ we have the expression
$$
\lambda_G(z)=\frac{|h'(z)|}{2\Im h(z)}
=\frac{|w\cdot dw/dz|}{2\Re w \Im w}.
$$
Note here that $dw/dz=(\alpha-\bar\alpha)/(z-\bar\alpha)^2
=2ir/(z-\bar\alpha)^2.$
Since
$$
w
=\frac{(z-\alpha)(\bar z-\alpha)}{|z-\bar\alpha|^2}
=\frac{U+iV}{|z-\bar\alpha|^2},
$$
where
$$
U=|z|^2-2\,\Re\alpha\cdot\Re z+\Re(\alpha^2)
=|z|^2+2c\,\Re z+c^2-r^2=|z+c|^2-r^2
$$
and
$$
V=-2\Im\alpha\cdot\Re z+\Im(\alpha^2)
=-2r\,\Re(z+c),
$$
Hence, we compute
\begin{align*}
\lambda_G(z)&
=\frac{2r|U+iV|/|z-\bar\alpha|^4}{2UV/|z-\bar\alpha|^4}
=\frac{r\sqrt{U^2+V^2}}{UV}.
\end{align*}
Recall that $a:=-c-r=\hat m(G)\in\partial\Omega.$
For a small enough $t>0,$ we choose $z=a-t=-c-r-t.$
Then, $U=(r+t)^2-r^2=t(2r+t),~V=2r(r+t)$ and thus
$$
\lambda_G(a-t)=\frac{(r+t)^2+r^2}{2t(2r+t)(r+t)}.
$$
Finally, we observe
\begin{align*}
\tau(a-t,a)\mu_G(a-t)
&=\frac{1+(c+r+t)^2}{1+(c+r)(c+r+t)}\cdot
\frac{(r+t)^2+r^2}{2(2r+t)(r+t)} \\
&=\frac12-\frac {1+c^2-r^2}{4r[1+(c+r)^2]}t+O(t^2).
\end{align*}
By \eqref{eq:cond}, we see that $\tau(a,a-t)\mu_G(a-t)<1/2$
for a small enough $t>0.$
Since $\mu_\Omega(z)\le\mu_G(z)$ for $z\in G,$ we observe that
$$
\Ch_\tau(\Omega)\le \tau(a-t,\partial\Omega)\mu_\Omega(a-t)
\le\tau(a-t,a)\mu_G(a-t)<1/2
$$
for such a $t>0.$
Hence we conclude that $\Ch_\tau(\Omega)<1/2$ as required.

\bigskip

{\bf Acknowledgements.}
The author is grateful to Professor Hajime Tanaka for useful conversations on the graph theory.
He would also like to thank the referee for his or her careful reading.

\def\cprime{$'$} \def\cprime{$'$} \def\cprime{$'$}
\providecommand{\bysame}{\leavevmode\hbox to3em{\hrulefill}\thinspace}
\providecommand{\MR}{\relax\ifhmode\unskip\space\fi MR }
\providecommand{\MRhref}[2]{%
  \href{http://www.ams.org/mathscinet-getitem?mr=#1}{#2}
}
\providecommand{\href}[2]{#2}


\begin{thebibliography}{10}

\bibitem{AW:SP}
F.~G. Avkhadiev and K.-J. Wirths, \emph{Schwarz-{P}ick {T}ype {I}nequalities},
  Frontiers in Mathematics, Birkh\"auser Verlag, Basel, 2009.

\bibitem{BP78}
A.~F. Beardon and {Ch}. Pommerenke, \emph{The {P}oincar\'e metric of plane
  domains}, J. London Math. Soc. (2) \textbf{18} (1978), 475--483.

\bibitem{Diestel:gt}
R.~Diestel, \emph{Graph {T}heory}, fourth ed., Graduate Texts in Mathematics,
  vol. 173, Springer, Heidelberg, 2010.

\bibitem{Duren:univ}
P.~L. Duren, \emph{Univalent {F}unctions}, Springer-Verlag, 1983.

\bibitem{HM92}
R.~Harmelin and D.~Minda, \emph{Quasi-invariant domain constants}, Israel J.
  Math. \textbf{77} (1992), 115--127.

\bibitem{KL:hg}
L.~Keen and N.~Lakic, \emph{Hyperbolic {G}eometry from a {L}ocal {V}iewpoint},
  Cambridge University Press, Cambridge, 2007.

\bibitem{Keogh76}
F.~R. Keogh, \emph{A characterization of convex domains in the plane}, Bull.
  London Math. Soc. \textbf{8} (1976), 183--185.

\bibitem{KM01}
S.~Kim and D.~Minda, \emph{The hyperbolic metric and spherically convex
  regions}, J. Math. Kyoto Univ. \textbf{41} (2001), 297--314.

\bibitem{MM90}
D.~Mej\'{\i}a and D.~Minda, \emph{Hyperbolic geometry in $ k $-convex regions},
  Pacific J. Math. \textbf{141} (1990), 333--354.

\bibitem{Minda85a}
D.~Minda, \emph{Estimates for the hyperbolic metric}, Kodai Math. J. \textbf{8}
  (1985), 249--258.

\bibitem{Minda86}
\bysame, \emph{The hyperbolic metric and {B}loch constants for spherically
  convex regions}, Complex Var. \textbf{5} (1986), 127--140.

\bibitem{Sugawa13}
T.~Sugawa, \emph{Spherical density of hyperbolic metric and uniform
  perfectness}, J. Analysis \textbf{21} (2013), 157--166.

\end{thebibliography}
\end{document}